\documentstyle[11pt]{article}
\title{Injectives in Residuated Algebras}
\author{{\sc Hector Freytes, Buenos Aires} \thanks{During the
preparation of this paper the author was supported by a
Fellowship from the FOMEC Program. The author expresses his
gratitude to  Roberto Cignoli, for his advice during the
preparation of this paper, and to the Referee for her/his many
suggestions to improve the presentation of this paper.} }
\date{}
\begin{document}
\bibliographystyle{plain}
\maketitle

\begin{abstract}
\noindent \noindent Injectives in several classes of structures
associated with logic are characterized. Among the classes
considered are residuated lattices, MTL-algebras, IMTL-algebras,
BL-algebras, NM-algebras and bounded hoops.
\end{abstract}
\begin{small}



\end{small}
\bibliography{pom}

\newtheorem{theo}{Theorem}[section]
\newtheorem{definition}[theo]{Definition}
\newtheorem{lem}[theo]{Lemma}
\newtheorem{prop}[theo]{Proposition}
\newtheorem{coro}[theo]{Corollary}
\newtheorem{exam}[theo]{Example}
\newtheorem{rema}[theo]{Remark}{\hspace*{4mm}}
\newtheorem{example}[theo]{Example}

\newcommand{\proof}{\noindent {\em Proof:\/}{\hspace*{4mm}}}
\newcommand{\qed}{\hfill$\Box$}
\newcommand{\ninv}{\mathord{\sim}} %involutive negation

\section*{Introduction} Residuated structures, rooted in the work
of Dedekind on the ideal theory of rings, arise in many fields of
mathematics, and are particularly common among algebras
associated with logical systems. They are  structures $\langle A,
\odot, \rightarrow, \leq \rangle$ such that $A$ is a nonempty
set, $\leq$ is a partial order on $A$ and $\odot$ and
$\rightarrow$ are binary operations such that the following
relation holds for each $a, b , c$ in $A$: $$a \odot b \leq c
\,\,\,\,\,  \mbox{{\em iff}} \,\,\,\,\,  a \leq b \rightarrow
c.$$ Important examples  of residuated structures related to
logic are Boolean algebras (corresponding to classical logic),
Heyting algebras (corresponding to intuitionism), residuated
lattices (corresponding to logics without contraction rule
\cite{TK}), BL-algebras (corresponding to H\'ajek's basic fuzzy
logic \cite{HAJ}), MV-algebras (corresponding to \L ukasiewicz
many-valued logic \cite{CDM}). All these examples, with the
exception of residuated lattices are {\em hoops\/}
\cite{BlokFer}, i. e., they satisfy the equation $x \odot (x
\rightarrow y) = y \odot (y \rightarrow x)$.

The aim of this paper is to investigate injectives and absolute
retracts in classes of residuated lattices and bounded hoops. In
\S2 and \S3 we also present some results on injectives in more
general varieties.

The paper is structured as follows.   In \S1 we recall some basic
definitions and properties.  In \S2 we show that under some mild
hypothesis on a variety $\cal V$ of algebras, the existence of
nontrivial injectives is equivalent to the existence of a
self-injective  maximum simple algebra. In \S3 we use ultrapowers
to obtain lattice properties of the injectives in varieties of
ordered algebras. The results of \S2 and \S3 are applied in \S4,
\S7 and in \S14 to the study of injectives in varieties of
residuated lattices, prelinear residuated lattices and bounded
hoops, respectively. In the remaining sections we consider
injectives in several subvarieties of residuated lattices which
appear in the literature.  The results obtained are summarized in
Table 1.

\section{Basic Notion}
We recall from from \cite{BD} and \cite{Bur}  some basic notion
of injectives and universal algebra. Let ${\cal A}$ be a class of
algebras. For all algebras $A, B$ in ${\cal A}$, $[A,B]_{\cal A}$
will denote the set of all homomorphism $g:A\rightarrow B$.  An
algebra $A$ in ${\cal A}$ is {\bf injective} iff for every
monomorphism  $f \in [B, C]_{{\cal A}}$ and every  $g \in [B,
A]_{{\cal A}}$ there  exists  $h \in [C, A]_{{\cal A}}$ such that
$hf = g$; $A$ is {\bf self-injective} iff every homomorphism from
a subalgebra of $A$ into $A$, extends to an endomorphism of $A$.
An algebra $B$ is a {\bf retract} of an algebra $A$ iff there
exists $g \in [B,  A]_{{\cal A}}$ and $f\in [A, B]_{{\cal A}}$
such that  $fg = 1_B$.  It is well known that a retract of an
injective object is injective. An algebra $B$ is called an {\bf
absolute retract} in ${\cal A}$ iff it is a retract of each of
its extensions in ${\cal A}$. For each algebra $A$, we denote by
$Con(A)$, the congruence lattice of $A$, the diagonal congruence
is denoted by $\Delta$ and the largest congruence $A^2$ is
denoted by $\nabla$. A congruence $\theta_M$ is said to be
maximal iff $\theta_M \not= \nabla$ and there is no congruence
$\theta$ such that $\theta_M \subset \theta \subset \nabla$. An
algebra $I$ is {\bf simple} iff $Con(I) = \{\Delta,\nabla \}$. A
nontrivial algebra $T$ is said to be {\bf minimal} in $\cal A$
iff for each nontrivial algebra $A$ in $\cal A$, there exists a
monomorphism $f:T\rightarrow A$. A simple algebra $I_M$ is said
to be {\bf maximum simple} iff for each simple algebra $I$, $I$
can be embedded in $I_M$. A simple algebra is {\bf hereditarily
simple} iff all its subalgebras are simple. An algebra $A$ is {\bf
semisimple} iff it is a subdirect product of simple algebras. An
algebra $A$ is {\bf rigid} iff the identity homomorphism is the
only automorphism. An algebra $A$ has the {\bf congruence
extension property} (CEP) iff for each subalgebra $B$ and $\theta
\in Con(B)$ there is a $\phi \in Con(A)$ such that $\theta = \phi
\cap A^2$.  A variety ${\cal V}$ satisfies CEP iff every algebra
in ${\cal V}$ has the CEP. It is clear that if ${\cal V}$
satisfies CEP then every simple algebra is hereditarily simple.

\section{Injectives and simple algebras}

\begin{definition}
{\rm Let $\cal V$ be a variety. Two constant terms  $0,1$ of the
language of $\cal V$ are called  {\it distinguished constants}
iff $A\models 0\not=1$ for each nontrivial algebra $A$ in $\cal
V$.}
\end{definition}

\begin{lem}\label{Disting}
Let $\cal A$ be variety with distinguished constants $0,1$ and let
$A$ be a nontrivial algebra in $\cal A$. Then $A$ has maximal
congruences, and for each simple algebra $I\in \cal A$, all
homomorphisms $f:I\rightarrow A$ are monomorphisms.

\end{lem}

\begin{proof}
Since for each homomorphism $f:A\rightarrow B$ such that $B$ is a
nontrivial algebra, $f(0)\not=f(1)$ then for each $\theta \in
Con(A)\backslash \{A^2\}$, $(1,0)\notin \theta$. Thus a standard
application of Zorn lemma shows that $Con(A)\backslash \{A^2\}$
has maximal elements. The second claim follows from the
simplicity of $I$ and $f(0)\not=f(1)$. \qed
\end{proof}

\begin{theo} \label{Injective simple}
Let ${\cal A}$ be a variety with distinguished constants $0, 1$
having a minimal algebra. If ${\cal A}$ has nontrivial
injectives, then there exists a maximum simple algebra $I$.
\end{theo}

\begin{proof}
Let $A$ be a nontrivial injective in ${\cal A}$. By
Lemma~\ref{Disting} there is a maximal congruence $\theta$ of $A$.
Let  $I = A/\theta$  and $p:A\rightarrow I$ be the canonical
projection. Since ${\cal A}$ has a minimal algebra it is clear
that for each simple algebra $J$, there exists a monomorphism
$h:J \rightarrow A$.  Then the composition $ph$ is a monomorphism
from $J$ into $I$. Thus $I$ is a maximum simple algebra. \qed
\end{proof}

\noindent
\\
We want to establish a kind of the converse of the above theorem.
\\

\begin{theo}\label{Simple injective}
Let ${\cal A}$ be a variety satisfying CEP, with distinguished
constants $0, 1$. If $I$ is a self-injective maximum simple
algebra in ${\cal A}$ then $I$ is injective.
\end{theo}

\begin{proof}
For each monomorphism $g:A\rightarrow B$ we consider the
following diagram in ${\cal A}$:

\begin{center}
\unitlength=1mm
\begin{picture}(60,20)(0,0)
\put(8,16){\vector(3,0){5}} \put(2,10){\vector(0,-2){5}}

\put(2,16){\makebox(0,0){$A$}} \put(20,16){\makebox(0,0){$I$}}
\put(2,0){\makebox(0,0){$B$}}

\put(2,20){\makebox(17,0){$f$}} \put(2,8){\makebox(-5,0){$g$}}
\end{picture}
\end{center}

By CEP, $I$ is hereditarily simple. Hence $f(A)$ is simple and
$Ker(f)$ is a maximal congruence of $A$ such that $(0,1) \notin
Ker(f)$. Further $Ker(f)$ can be extended to a maximal congruence
$\theta$ in $B$. It is clear that  $(0,1)\notin \theta$ and
$\theta \cap A^2 = Ker(f)$. Thus if we consider the canonical
projection $p:B \rightarrow B/\theta$, then there exists a
monomorphism $g':f(A) \rightarrow B/\theta$ such that

\begin{center}
\unitlength=1mm
\begin{picture}(60,20)(0,0)
\put(8,16){\vector(3,0){5}} \put(2,10){\vector(0,-2){5}}
\put(8,0){\vector(1,0){5}} \put(20,10){\vector(0,-2){5}}

\put(26,16){\vector(3,0){5}}

\put(2,16){\makebox(0,0){$A$}} \put(20,16){\makebox(0,0){$f(A)$}}
\put(2,0){\makebox(0,0){$B$}}
\put(20,0){\makebox(0,0){$B/\theta$}}

\put(36,16){\makebox(0,0){$I$}}

\put(1,10){\makebox(17,-3){$\equiv$}}

\put(2,20){\makebox(17,0){$f$}} \put(2,8){\makebox(-5,0){$g$}}
\put(14,-4){\makebox(-5,2){$p$}} \put(19,7){\makebox(-5,2){$g'$}}

\put(20,20){\makebox(17,0){$1_{f(A)}$}}
\end{picture}
\end{center}

Since $I$ is maximum simple, $B/\theta$ is isomorphic to a
subalgebra of $I$. Therefore, since that $I$ is self-injective,
there exists a monomorphism $\varphi: B/\theta \rightarrow I$
such that $\varphi g' = 1_f(A)$. Thus $(\varphi p)g = f$ and $I$
is injective. \qed
\end{proof}

\begin{lem}\label{RIGID}
 If $A$ is a rigid simple injective
 algebra in a variety,  then all the subalgebras of $A$ are rigid. \qed
\end{lem}

\section{Injectives, ultrapowers and lattice properties}

\noindent We recall from \cite{Bir} some basic notions on ordered
sets that will play an important role in what follows. An ordered
set $L$ is called {\bf bounded} provided it has  a smallest
element $0$ and a greatest element $1$. The {\bf decreasing
segment} ${\rm ({\it a}]}$ of $L$ is defined as the set $\{\,x\in
L : x\leq a \}$. The increasing segment $[a)$ is defined dualy. A
subset $X$ of $L$ is called {\bf down directed} ({\bf upper
directed}) iff for all $a,b\in X$, there exists $x\in X$ such that
$x \leq a$ and $x \leq b$ ($a \leq x$ and $b \leq x$).\\

\begin{lem}
Let $L$ be a lattice and $X$ be a down (upper) directed subset of
$L$ such that $X$ does not have a minimum (maximum) element. If
$\cal F$ is the filter in ${\cal P}(X)$ generated by the
decreasing (increasing) segments of $X$, then there exists a
nonprincipal ultrafilter $\cal U$ such that $\cal F \subseteq \cal
U$.
\end{lem}

\begin{proof}
Let ${\rm ({\it a}]}$, ${\rm ({\it b}]}$ be decreasing segments
of $X$. Since $X$ is a down directed subset, there exists $x\in
X$ such that $x\leq a$ and $x\leq b$, whence $x\in ({\it a}]\cap
({\it b}]$ and $\cal F$ is a proper filter of ${\cal P}(X)$. By
the ultrafilter theorem there exists an ultrafilter $\cal U$ such
that $\cal F \subseteq \cal U$. Suppose that $\cal U$ is the
principal filter generated by ${\rm ({\it c}]}$. Since $X$ does
not have a minimum element, there exists $x\in X$ such that
$x<c$. Thus ${\rm ({\it x}]} \in \cal U$ and it is a proper subset
of ${\rm ({\it c}]}$, a contradiction. Hence $\cal U$ is not a
principal filter. By duality, we can establish the same result
when $X$ is an upper directed set. \qed
\end{proof}

\begin{definition}\label{LA}
{\rm A variety ${\cal V}$ of algebras has {\it lattice-terms} iff
there are terms of the language of ${\cal V}$ defining on each
$A\in {\cal V}$ operations $\lor$, $\land$, such that $\langle
A,\lor,\land\rangle$ is a lattice. ${\cal V}$ has {\it bounded
lattice-terms} if, moreover, there are two constant terms $0$,$1$
of the language of ${\cal V}$ defining on each $A\in {\cal V}$ a
bounded lattice $\langle A,\lor,\land, 0, 1\rangle$. The order in
$A$, denoted by $L(A)$, is called the {\it natural order} of $A$.}
\end{definition}

\noindent Observe that each subvariety of a variety with
(bounded) lattice-terms is also a variety with (bounded)
lattice-terms.
\\

\noindent Let ${\cal V}$ be a variety with lattice-terms and
$A\in {\cal V}$. ${A^X}/{\cal U}$ will always denote the
ultrapower corresponding to a down (upper) directed set $X$ of
$A$ with respect to the natural order, without smallest
(greatest) element and a nonprincipal ultrafilter ${\cal U}$ of
${\cal P}(X)$, containing the filter generated by the decreasing
(increasing) segments of $X$. For each $f\in A^X$, $[f]$ will
denote the ${\cal U}$-equivalence class of $f$. Thus $[1_X]$ is
the ${\cal U}$-equivalence class of the canonical injection
$X\hookrightarrow A$ and for each $a\in A$, $[a]$ is  the ${\cal
U}$-equivalence class of the constant function $a$ in $A^X$. It
is well known that $i_A(a) = [a]$ defines a monomorphism
$A\rightarrow {A^X}/{\cal U}$  {\rm (see \cite[Corollary
4.1.13]{CK})}.

\begin{theo}\label{Ultraproducto}
Let ${\cal V}$ be a variety with lattice-terms. If there exists
an absolute retract $A$ in ${\cal V}$, then each down directed
subset $X\subseteq A$ has an infimum, denoted by $\bigwedge X$.
Moreover if $P(x)$ is a first-order positive formula {\rm (see
\cite{CK})} of the language of ${\cal V}$ such that each $a\in X$
satisfies $P(x)$, then $\bigwedge X$ also satisfies $P(x)$.
\end{theo}

\begin{proof}
Let $X$ be a down directed subset of the absolute retract $A$.
Suppose that $X$ does not admit a minimum element and consider an
ultrapower ${A^X}/{\cal U}$. Since $A$ is an absolute retract
there exists a homomorphism $\varphi$ such that the following
diagram is commutative:

\begin{center}
\unitlength=1mm
\begin{picture}(60,20)(0,0)
\put(8,16){\vector(3,0){5}} \put(2,10){\vector(0,-2){5}}
\put(8,2){\vector(1,1){7}} \put(2,16){\makebox(0,0){$A$}}
\put(20,16){\makebox(0,0){$A$}}
\put(2,0){\makebox(0,0){${A^X}/{\cal U}$}}
\put(2,20){\makebox(17,0){$1_A$}}
\put(2,10){\makebox(13,0){$\equiv$}}
\put(2,8){\makebox(-5,0){$i_A$}}
\put(16,0){\makebox(-5,2){$\varphi$}}
\end{picture}
\end{center}

We first prove that $\varphi ([1_X])$ is a lower bound of $X$.
Let $a\in X$. Then $[1_X]\leq [a]$ since $\{x\in X : 1_X(x)\leq
a(x)   \} = \{x\in X : x\leq a\} \in \cal U$. Thus $\varphi
([1_X])\leq \varphi([a]) = a$ and $\varphi ([1_X])$ is a lower
bound of $X$. We proceed now to prove that $\varphi ([1_X])$ is
the greatest lower bound of $X$. In fact, if $b \in A$ is a lower
bound of $X$ then for each $x\in X$ we have $b\leq x$. Thus
$[b]\leq [1_X]$ since $\{x\in X : b(x)\leq 1_X(x)\} = \{x\in X :
b\leq x\} = X \in \cal U$. Now we have $b = \varphi ([b])\leq
\varphi ([1_X])$. This proves that $\varphi ([1_X])= \bigwedge
X$. If each $a\in X$ satisfies the first order formula $P(x)$
then $[1_X]$ satisfies $P(x)$ and, since $P(x)$ is a positive
formula, it follows from {\rm (\cite[Theorem 3.2.4]{CK} )} that
$\varphi ([1_X])$ satisfies $P(x)$. \qed
\end{proof}
\\

\noindent In the same way, we can establish the dual version of
the above theorem. Recalling that a lattice is complete iff there
exists the infimum $\bigwedge X$ (supremum $\bigvee X$), for each
down directed (upper directed) subset $X$, we have the following
corollary:

\begin{coro}
Let ${\cal V}$ be a variety with lattice-terms. If $A$ is an
absolute retract in ${\cal V}$, then $L(A)$ is a complete
lattice. \qed
\end{coro}

\section{Residuated Lattices and Semisimplicity}

\begin{definition}
{\rm A  {\it residuated lattice} {\rm \cite{TK}} or {\it
commutative integral residuated $0,1$-lattice}  {\rm \cite{JT}},
is an algebra $ \langle A, \land, \lor, \odot, \rightarrow, 0, 1
\rangle$ of type $ \langle 2, 2, 2, 2, 0, 0 \rangle$ satisfying
the following axioms:

\begin{enumerate}
\item
$\langle A,\odot,1 \rangle$ is an abelian monoid,

\item
$L(A) = \langle A, \lor, \land, 0,1 \rangle$ is a bounded lattice,

\item
$(x \odot y)\rightarrow z = x\rightarrow (y\rightarrow z)$,

\item
$((x\rightarrow y)\odot x)\land y = (x\rightarrow y)\odot x$,

\item
$(x\land y)\rightarrow y = 1$.

\end{enumerate}

\noindent $A$ is  called an {\it involutive residuated lattice} or
{\it Girard monoid} {\rm \cite{UHO}} if it also satisfies the
equation:

\begin{enumerate}
\item[6.]
$(x\rightarrow 0)\rightarrow 0 = x$.
\end{enumerate}

\noindent $A$ is called {\it distributive} if  satisfies 1. -- 5.
as well as:

\begin{enumerate}
\item[7.]
$x\land (y \lor z) = (x\land y)\lor (x\land z)$.
\end{enumerate}

}
\end{definition}

\noindent The variety of residuated lattices is denoted by ${\cal
RL}$, and the subvariety of Girad monoids is noted by ${\cal GM}$.
Following the notation used in {\rm \cite{JT}}, the variety of
residuated lattices that satisfy the distributive law is denoted
by ${\cal DRL}$, and ${\cal DGM}$ will denote the variety of
distributive Girad monoids. It is clear that $0,1$ are
distinguished constant terms in ${\cal RL}$. Moreover, $\{0,1\}$
is a subalgebra of each nontrivial $A \in {\cal RL}$, which is a
boolean algebra. Hence $\{0,1\}$ with its natural boolean algebra
structure is the minimal algebra in each nontrivial subvariety of
${\cal RL}$. Thus the variety ${\cal BA}$ of boolean algebras is
contained in all nontrivial varieties of residuated lattices. On
each residuated lattice $A$ we can define a unary operation
$\neg$ by $\neg x = x\rightarrow 0$.  We also define for all
$a\in A$, $a^1 = a$ and $a^{n+1}= a^n\odot a$. An element $a$ in
$A$ is called {\bf idempotent} iff $a^2 = a$, and it is called
{\bf nilpotent} iff there exists a natural number $n$ such that
$a^n= 0$. The minimum $n$ such that $a^n= 0$ is called {\bf
nilpotence order} of $a$. An element $a$ in $A$ is called {\bf
dense} iff $\neg a = 0$ and it is called a {\bf unity} iff for
all natural numbers $n$, $\neg (a^n)$ is nilpotent. The set of
dense elements of $A$ will be denoted by $Ds(A)$. We recall now
some well-known facts about implicative filters and congruences
on residuated lattices. Let A be a residuated lattice and
$F\subseteq A$. Then $F$ is an {\bf implicative filter} iff it
satisfies the following conditions:

\begin{enumerate}
\item
$1\in F$,

\item
if $x\in F$ and $x\rightarrow y \in F$ then $y\in F$.

\end{enumerate}

\noindent It is easy to verify that a nonempty subset $F$ of a
residuated lattice $A$ is an implicative filter iff for all
$a,b\in A$:

\begin{enumerate}
\item[-]
If $a\in F$ and $a\leq b$ then $b\in F$,

\item[-]
if $a,b\in F$ then $a\odot b \in F$.

\end{enumerate}

Note that an implicative filter $F$ is proper iff $0$ does not
belong to $F$. The intersection of any family of implicative
filters of $A$ is again an implicative filter of $A$. We denote
by $ \langle X \rangle$ the implicative filter generated by
$X\subseteq A$, i.e., the intersection of all implicative filters
of $A$ containing $X$. We abbreviate this as $ \langle a \rangle$
when $X=\{a\}$ and it is easy to verify that $\langle X \rangle =
\{x\in A: \exists \hspace{0.2 cm} w_1 \cdots w_n \in X
\hspace{0.2 cm}  \mbox{such that}  \hspace{0.2 cm} x\geq
w_1,\odot \cdots , \odot w_n \}$. For any implicative filter $F$
of $A$, $\theta_F = \{(x,y)\in A^2: x\rightarrow y, y\rightarrow
x \in F \}$ is a congruence on $A$.

\noindent Moreover $F = \{ x \in A : (x,1)\in \theta_F \}$.
Conversely, if $\theta \in Con(A)$ then $F_{\theta} = \{x \in A :
(x,1) \in \theta \}$ is an implicative filter and $(x,y) \in
\theta $ iff $(x \rightarrow y, 1) \in \theta$ and $(y
\rightarrow x, 1) \in \theta$. Thus the correspondence $F
\rightarrow \theta_F$ is a bijection from the set of implicative
filters of $A$ onto the set $Con(A)$. If $F$ is an implicative
filter of $A$, we shall write $A/F$ instead of $A/\theta_F$, and
for each $x\in A$ we shall write  $x/\theta_F$ for the equivalence
class of $x$.

\begin{prop}
If ${\cal A}$ is a subvariety of ${\cal RL}$, then ${\cal A}$
satisfies CEP.
\end{prop}

\begin{proof}
This follows from the same argument used in {\rm (\cite[Theorem
1.8]{BlokFer})}.
\qed\\
\end{proof}

\noindent If $A$ is a residuated lattice then we define

$$Rad(A) = \bigcap \{F: F \hspace{0.2 cm} is \hspace{0.2 cm} a \hspace{0.2 cm} maximal\hspace{0.2 cm} implicative\hspace{0.2 cm} filter\hspace{0.2 cm}  in \hspace{0.2 cm} A\}.$$

\noindent It is clear that $A$ is semisimple iff $Rad(A) =
\{1\}$. If ${\cal A}$ is a subvariety of ${\cal RL}$, we denote
by ${\cal S}em({\cal A}) $ the subclass of ${\cal A}$ whose
elements are the semisimple algebras of ${\cal A}$. Thus we have
${\cal S}em({\cal A}) = \{A/Rad(A):A\in {\cal A} \}$.

\begin{prop}\label{Rad}
Let A be a residuated lattice. Then:

\begin{enumerate}
\item
$A$ is simple iff for each $a<1$, $a$ is nilpotent.

\item
$Rad(A) = \{a\in A : a \hspace{0.2 cm} is \hspace{0.2 cm} unity\}
$.

\item
$Ds(A)$ is an implicative filter in $A$ and $Ds(A) \subseteq
Rad(A)$ .

\end{enumerate}
\end{prop}

\begin{proof}
1)\hspace{0.2 cm} Trivial. 2)\hspace{0.2 cm} See {\rm
(\cite[Lemma 4.6]{UHO} )}. 3) \hspace{0.2 cm}Follows immediately
from 2. \qed
\\
\end{proof}

\noindent If $Rad(A)$ has a least element $a$, i.e., $Rad(A)=
[a)$, then $a$ is called the {\bf principal unity} of $A$. It is
clear that a principal unity is an idempotent element and that it
generates the radical.

\begin{lem}\label{IMPUNIT}
Let $A$ be a residuated lattice having a principal unity $a$. If
$x\in Rad(A)$, then $x\rightarrow \neg a = \neg a$.
\end{lem}

\begin{proof}
$x\rightarrow \neg a = \neg (x\odot a) =  \neg a$ since $a$ is
the minimum unity. \qed
\end{proof}

\begin{prop}\label{NEGUNIT}
Let $A$ be  a linearly ordered residuated lattice. Then:

\begin{enumerate}
\item
$a$ is a unity in $A$ iff \hspace{0.1 cm} $a$ is not a nilpotent
element.

\item
If $a$ is a unity in $A$, then $\neg a < a$.

\end{enumerate}
\end{prop}

\begin{proof}
1)\hspace{0.2 cm} If $a < 1$ and there exists a natural number
$n$ such that $a^n = 0$, then $\neg(a^n) = 1$ and $a$ is not a
unity. Conversely, suppose $a$ is not a unity. Since $A$ is
linearly ordered, we must have $a^n \leq \neg \neg(a^n) <
\neg(a^n)$. Hence $a^{2n} = 0$ and $a$ is nilpotent, which is a
contradiction. 2)\hspace{0.2 cm} Is an obvious consequence of 1).
\qed
\end{proof}

\begin{coro}\label{UNMTL}
Let $A$ be a residuated lattice such that there exists an
embedding $f:A \rightarrow \prod_{i\in I} L_i$, with $L_i$ a
linearly ordered residuated lattice for each $i\in I $. Then $a$
is a unity in $A$ iff for each $i\in I$, $a_i =\pi_i f(a)$ is a
unity in $L_i$, where $\pi_i$ is the $i{\rm th}$-projection onto
$L_i$ .
\end{coro}

\begin{proof}
If $a$ is a unity in $A$ then $a_i = \pi_i f(a)$ is a unity in
$L_i$, because homomorphisms preserve unities. Conversely,
suppose that $a$ is not a unity. Therefore there is an $n$ such
that $\neg (a^n)$ is not nilpotent, and hence $\neg (a^n) \not
\leq \neg \neg (a^n)$.  Since $f$ is an embeding and since $L_i$
is linearly ordered for each $i\in I$, there exists $j\in I$ such
that $\neg\neg(a_j^n) \leq \neg (a_j^n)$, and by Proposition
\ref{NEGUNIT} $a_j$ is not a unity in $L_j$. \qed
\end{proof}

\begin{prop}\label{RadSub}
Let ${\cal A}$ be a subvariety of ${\cal RL}$. Then ${\cal
S}em({\cal A})$ is a reflective  subcategory, and the reflector
{\rm \cite{BD}} preserves monomorphism.
\end{prop}

\begin{proof}
If $A \in {\cal A}$, for each $x\in A$, $[x]$ will denote the
$Rad(A)$-congruence class of $x$. We define ${\cal S}(A)=
A/Rad(A)$, and for each $f \in [A,A']_{\cal A}$, we let ${\cal
S}(f)$ be defined by ${\cal S}(f)([x]) = [f(x)]$ for each $x\in
A$. Since homomorphisms preserve unity, we obtain a well defined
function ${\cal S}(f): A/Rad(A)\rightarrow A'/Rad(A')$. It is easy
to check that ${\cal S}$ is a functor from ${\cal A}$ to ${\cal
S}em({\cal A})$. To show that ${\cal S}$ is a reflector, note
first that if $p_A:A \rightarrow A/Rad(A)$ is the canonical
projection, then the following diagram is commutative:

\begin{center}
\unitlength=1mm
\begin{picture}(20,20)(0,0)
\put(8,16){\vector(3,0){5}} \put(-4,10){\vector(0,-2){5}}
\put(8,0){\vector(1,0){5}} \put(24,10){\vector(0,-2){5}}

\put(-4,16){\makebox(0,0){$A$}} \put(24,16){\makebox(0,0){$A'$}}
\put(-6,0){\makebox(0,0){$A/Rad(A)$}}
\put(27,0){\makebox(0,0){$A/Rad(A')$}}

\put(2,10){\makebox(17,-3){$\equiv$}}

\put(2,20){\makebox(17,0){$f$}} \put(-7,8){\makebox(-5,0){$p_A$}}
\put(14,-5){\makebox(-5,2){${\cal S}(f)$}}
\put(33,8){\makebox(-5,2){$p_{A'}$}}
\end{picture}
\end{center}

\hspace{0.2 cm}

Suppose that $B\in {\cal S(A)}$ and $f\in [A,B]_{{\cal A}}$.
Since $Rad(B) = \{1\}$, the mapping $[x]\mapsto f(x)$  defines a
homomorphism $g:A/Rad(A)\rightarrow B$ that makes  the following
diagram  commutative:

\begin{center}
\unitlength=1mm
\begin{picture}(20,20)(0,0)
\put(8,16){\vector(3,0){5}} \put(2,10){\vector(0,-2){5}}
\put(10,4){\vector(1,1){7}}

\put(2,10){\makebox(13,0){$\equiv$}}

\put(2,16){\makebox(0,0){$A$}} \put(20,16){\makebox(0,0){$B$}}
\put(2,0){\makebox(0,0){$A/Rad(A)$}}
\put(2,20){\makebox(17,0){$f$}} \put(2,8){\makebox(-6,0){$p_A$}}
\put(18,2){\makebox(-4,3){$g$}}
\end{picture}
\end{center}

\noindent and it is obvious that $g$ is the only homomorphism  in
$[A/Rad(A),B]_{{\cal S}em({\cal A})}$ making the triangle
commutative. Therefore we have proved that ${\cal S}$ is a
reflector. We proceed to prove that ${\cal S}$ preserves
monomorphisms. Let $f\in [A,B]_{{\cal A}}$ be a monomorphism and
suppose that $({\cal S}(f))(x)=({\cal S}(f))(y)$, i.e.,
$[f(x)]=[f(y)]$. Then for each number $n$ there exists a number
$m$ such that $0 = (\neg ((f(x)\rightarrow f(y))^n))^m = f((\neg
((x\rightarrow y)^n))^m)$. Since $f$ is a monomorphism then $(\neg
((x\rightarrow y)^n))^m = 0$ and $x\rightarrow y \in Rad(A)$.
Interchanging $x$ and $y$, we obtain $[x]= [y]$ and ${\cal S}(f)$
is a monomorphism. \qed
\end{proof}

\begin{coro}\label{INJBOT}
Let ${\cal A}$ be a subvariety of ${\cal RL}$. If $A$ is
injective in ${\cal S}em({\cal A})$ then $A$ is injective in
${\cal A}$.
\end{coro}

\begin{proof}
It is well known that if ${\cal D}$ is a reflective subcategory
of ${\cal A}$ such that the reflector  preserves monomorphisms
then an injective object in ${\cal D}$ is also injective in
${\cal A}$\hspace{0.2 cm}{\rm \cite[I.18]{BD}}. Then this theorem
follows from Propositions \ref{RadSub}.
\qed \\
\end{proof}

\noindent We will say that a variety ${\cal A}$ is ${\bf
radical-dense}$ provided that ${\cal A}$ is a subvariety of
${\cal RL}$ and $Rad(A) = Ds(A)$ for each $A$ in ${\cal A}$. An
example of a radical-dense variety is the variety ${\cal H}$ of
Heyting algebras (i.e., ${\cal RL}$ plus the equation $x\odot y =
x\land y$).

\begin{theo}\label{H3}
Let ${\cal A}$ be a radical-dense variety. If $A$ is a
non-semisimple absolute retract in ${\cal A}$, then $A$ has a
principal unity $\epsilon$ and  $\{0,\epsilon,1\}$ is a
subalgebra of $A$ isomorphic to the three element Heyting algebra
$H_3$.
\end{theo}

\begin{proof}
Let $A$ be a non-semisimple absolute retract. Unities are
characterized by the first order positive formula $\neg x = 0$
because $Rad(A)= Ds(A)$. Since $Ds(A)$ is a down-directed set, by
Theorem \ref{Ultraproducto} there exists a minimum dense element
$\epsilon$. It is clear that $\epsilon$ is the principal unity
and since $\epsilon < 1$,  $\{0,\epsilon,1\}$ is a subalgebra of
$A$, which coincides with the three element Heyting algebra $H_3$.
\qed
\end{proof}

\begin{definition}
{\rm Let ${\cal A}$ be a radical-dense variety. An algebra $T \in
{\cal A}$ is called a {\it $test_d$-algebra} iff there are
$\epsilon, t \in Rad(T) $ such that $\epsilon$ is an idempotent
element, $t < \epsilon$ and $\epsilon \rightarrow t \leq
\epsilon$}.
\end{definition}

\noindent An important example of a $test_d$-algebra is the
totally ordered four element Heyting algebra $H_4 = \{0 < b < a <
1\}$ whose operations are given as follows:

$$
x\odot y = x\land y,
$$

$$
x\rightarrow y = \cases {1, & if $x \leq y$,\cr y, & if $x > y.
$\cr }
$$

\begin{theo}\label{TEST1}
Let ${\cal A}$ be a radical-dense variety. If ${\cal A}$ has a
nontrivial injective and contains a $test_d$-algebra $T$, then
all injectives in ${\cal A}$ are semisimple.
\end{theo}

\begin{proof}
Suppose that there exists a non-semisimple injective $A$ in
${\cal A}$. Then by Lemma \ref{H3}, there is a monomorphism
$\alpha: H_3 \rightarrow A$ such that $\alpha(a)$ is the
principal unity in $A$. Let $i:H_3\rightarrow T$ be the
monomorphism such that $i(a)= \epsilon$. Since $A$ is injective,
there exists a homomorphism $\varphi:T\rightarrow A$ such that
the following diagram commutes

\begin{center}
\unitlength=1mm
\begin{picture}(20,20)(0,0)
\put(8,16){\vector(3,0){5}} \put(2,10){\vector(0,-2){5}}
\put(8,2){\vector(1,1){7}}

\put(2,10){\makebox(13,0){$\equiv$}}

\put(2,16){\makebox(0,0){$H_3$}} \put(20,16){\makebox(0,0){$A$}}
\put(2,0){\makebox(0,0){$T$}} \put(2,20){\makebox(17,0){$\alpha$}}
\put(2,8){\makebox(-6,0){$i$}}
\put(16,0){\makebox(-4,3){$\varphi$}}
\end{picture}
\end{center}

\noindent Since $\alpha(a)$ is the principal unity in $A$ and
$t\leq \epsilon$, then, by commutativity,  $\varphi(\epsilon) =
\varphi(t)= \alpha(a)$. Thus $\varphi(\epsilon \rightarrow t) =
1$, which is a contradiction since by hypothesis $\varphi(\epsilon
\rightarrow t) \leq \varphi(\epsilon) = \alpha(a) < 1$. Hence
${\cal A}$ has only semisimple injectives. \qed
\end{proof}

\section {Injectives in ${\cal RL}$, ${\cal GM}$, ${\cal DRL}$ and ${\cal DGM}$ }

\begin{prop}\label{EXT}
Let $A$ be a residuated lattice. Then the set  $A^\diamond =
\{(a,b)\in A\times A: a\leq b \}$ equipped with  the  operations \\

$(a_1,b_1)\land (a_2,b_2):= (a_1\land a_2, b_1\land b_2)$,

$(a_1,b_1)\lor (a_2,b_2):= (a_1\lor a_2, b_1\lor b_2)$,

$(a_1,b_1)\odot (a_2,b_2):= (a_1\odot a_2, (a_1\odot b_2)\lor
(a_2\odot b_1))$,

$(a_1,b_1)\rightarrow (a_2,b_2):= ((a_1\rightarrow a_2)\land (b_1\rightarrow b_2), a_1\rightarrow b_2)$.\\

\noindent is a residuated lattice, and the following properties
hold:

\begin{enumerate}

\item
The map $i:A\rightarrow A^\diamond$ defined by $i(a)=(a,a)$ is a
monomorphism.

\item
$\neg(a,b)= (\neg b, \neg a)$ and $\neg(0,1)=(0,1)$.

\item
$A$ is a Girard monoid iff $A^\diamond$ is a Girard monoid.

\item
$A$ is distributive iff $A^\diamond$ is distributive.

\end{enumerate}
\end{prop}

\begin{proof}
See {\rm \cite[IV Lemma 3.2.1]{UHO}}. \qed
\end{proof}

\begin{definition}
We say that a subvariety ${\cal A}$ of ${\cal RL}$ is {\bf
$\diamond$-closed} iff for all $A \in {\cal A}$, $A^\diamond \in
{\cal A}$.
\end{definition}

\begin{theo} \label{INJRL}
If a subvariety ${\cal A}$ of ${\cal RL}$ is {\bf
$\diamond$-closed},  then ${\cal A}$ has only trivial absolute
retracts.
\end{theo}

\begin{proof}
Suppose that there exists a non-trivial absolute retract $A$ in
${\cal A}$. Then by Proposition \ref{EXT} there exists an
epimorphism $f:A^\diamond \rightarrow A$ such that the following
diagram is commutative

\begin{center}
\unitlength=1mm
\begin{picture}(20,20)(0,0)
\put(8,16){\vector(3,0){5}} \put(2,10){\vector(0,-2){5}}
\put(8,2){\vector(1,1){7}}

\put(2,10){\makebox(13,0){$\equiv$}}

\put(2,16){\makebox(0,0){$A$}} \put(20,16){\makebox(0,0){$A$}}
\put(2,0){\makebox(0,0){$A^\diamond$}}
\put(2,20){\makebox(17,0){$1_A$}} \put(2,8){\makebox(-6,0){$i$}}
\put(16,0){\makebox(-4,3){$f$}}
\end{picture}
\end{center}

Thus there exists $a\in A$ such that $f(0,1)= a = f(a,a)$. Since
$(0,1)$ is a fixed point of the negation in $A^\diamond$ it
follows that $0<a<1$. We have $f(a,1)=1$. Indeed,
$(0,1)\rightarrow (a,a) = ((0 \rightarrow a)\land (1\rightarrow
a), 0 \rightarrow a) = (a,1)$. Thus $f(a,1)= f((0,1)\rightarrow
(a,a))= f(0,1)\rightarrow f(a,a)= a\rightarrow a =1$. In view of
this we have $1 = f(a,1)\odot f(a,1) = f((a,1) \odot(a,1)) =
f(a\odot a, (a\odot 1)\lor (a\odot 1))= f((a\odot a , a)) \leq
f((a,a))= a$, which is a contradiction since $a<1$. Hence ${\cal
A}$ has only trivial absolute retracts. \qed
\end{proof}

\begin{coro}
${\cal RL}$, ${\cal GM}$, ${\cal DRL}$ and ${\cal DGM}$  have
only trivial absolute retracts and injectives. \qed
\end{coro}

\section {Injectives in SRL-algebras}

\begin{definition}
A SRL-algebra is a residuated lattice satisfying the
equation:$$x\land \neg x = 0 \leqno(S)$$ The variety of
SRL-algebras is denoted by ${\cal SRL}$.
\end{definition}

\begin{prop}\label{NILPS}
If $A$ is a SRL-algebra, then $0$ is the only nilpotent in $A$.
\end{prop}

\begin{proof}
Suppose that there exists a nilpotent element $x$ in $A$ such
that $0<x$, having nilpotence order equal to $n$. By the
residuation property we have $x^{n-1}\leq \neg x$. Thus
\hspace{0.1cm} $x^{n-1} = x\land x^{n-1}\leq  x \land \neg x =
0$, which is a contradiction since $x$ has nilpotence order equal
to $n$. \qed
\end{proof}

\begin{coro}\label{SIMPSRL}
Let ${\cal A}$ be a subvariety of ${\cal SRL}$. Then the
two-element boolean algebra is the maximum simple algebra in
${\cal A}$ and ${\cal S}em({\cal A})  = {\cal BA}$.
\end{coro}

\begin{proof}
Follows from Propositions \ref{NILPS} and \ref{Rad}. \qed
\end{proof}

\begin{coro}\label{SRLRD}
If ${\cal A}$ is a subvariety of ${\cal SRL}$ then ${\cal A}$ is
a radical-dense variety.
\end{coro}

\begin{proof}
Let $A$ be an algebra in ${\cal A}$ and let $a$ be a unity. Thus
$\neg a$ is nilpotent and hence $\neg a = 0$. \qed
\end{proof}

\begin{coro}\label{SRLINJ}
If ${\cal A}$ is a subvariety of ${\cal SRL}$, then all complete
boolean algebras are injectives in ${\cal A}$.
\end{coro}

\begin{proof}
By Corollary \ref{SIMPSRL} the two-element boolean algebra is the
maximun simple algebra in ${\cal A}$. Since it is self-injective,
by Theorem \ref{Simple injective} it is injective. Since complete
boolean algebras are the retracts of powers of the two-element
boolean algebra, the result is proved. \qed
\end{proof}

\vspace{0.5cm}

\noindent As an application of this theorem we prove the
following results :

\begin{coro}
In ${\cal SRL}$ and ${\cal H}$, the only injectives  are complete
boolean algebras.
\end{coro}

\begin{proof}
Follows from Corollary \ref{SRLINJ} and Theorem \ref{TEST1}
because the $test_d$-algebra $H_4$ belongs to both varieties. \qed
\end{proof}

\begin{rema}
{\rm The fact that injective Heyting algebras are exactly
complete boolean algebras was proved in {\rm \cite{BH}} by
different arguments.}
\end{rema}

\section {MTL-algebras and absolute retracts}

\begin{definition}
{\rm An {\it MTL-algebra}  {\rm \cite{GE1}} is a residuated
lattice satisfying the pre-linearity equation
$$(x\rightarrow y)\lor (y\rightarrow x) = 1 \leqno(Pl)$$

\noindent The variety of MTL-algebras is denoted by ${\cal MTL}$}.
\end{definition}

\begin {prop}\label{SUBMTL}
Let $A$ be a residuated lattice. Then the following conditions are
equivalent:

\begin{enumerate}
\item
$A \in {\cal MTL}$.

\item
$A$ is a subdirect product of linearly ordered residuated
lattices.
\end{enumerate}
\end {prop}

\begin{proof}
{\rm \cite[Theorem 4.8 p. 76 ]{UHO}}. \qed
\\
\end{proof}

\begin {coro}\label{dist}
${\cal MTL}$ is subvariety of ${\cal DRL}$. \qed
\end {coro}

\begin {coro}\label{SIMMTL}
Let $A$ be a MTL-algebra.

\begin{enumerate}

\item
If $A$ is simple, then $A$ is linearly ordered.

\item
If $e$ is a unity in $A$, then $\neg e < e$.

\end{enumerate}
\end {coro}

\begin{proof}
1)\hspace{0.2 cm} Is an immediate consequence of Proposition
\ref{SUBMTL}. 2)\hspace{0.2 cm}If we consider that the {\it
i}{\rm th}-coordinate $\pi_if(e)$ of $e$ in the subdirect product
$f:A \rightarrow \prod_{i\in I} L_i$ is a unity, for each $i\in
I$, then by Proposition~\ref{NEGUNIT}, $\neg \pi_if(e) <
\pi_if(e)$. Thus $\neg e < e$. \qed
\\
\end{proof}

To obtain the analog of Theorem \ref{H3} for varieties of
MTL-algebras, we cannot use directly Theorem \ref{Ultraproducto},
because the property of being a unity is not a first order
property. We need to adapt the proof of Theorem \ref{H3} to this
case:

\begin {theo}\label{MINUNITY}
Let ${\cal A}$ be a subvariety of ${\cal MTL}$. If $A$ is an
absolute retract in ${\cal A}$ then $A$ has a principal unity $e$
in $A$.
\end {theo}

\begin{proof}
By Proposition \ref{SUBMTL} we can consider  a subdirect embedding
$f:A \rightarrow \prod_{i\in I} L_i$ such that $L_i$ is linearly
ordered. We define a family $H(L_i)$ in ${\cal A}$ as follows:
for each $i\in I$

\begin{enumerate}

\item[(a)]
if there exists $e_i = min\{u\in L_i : u \hspace{0.1 cm} is
\hspace{0.1 cm} unity\}$ then $H(L_i) = L_i$,

\item[(b)]
otherwise, $X = \{u\in L_i : u \hspace{0.1 cm} is \hspace{0.1 cm}
unity\}$ is a down-directed set without least element. Then by
Proposition \ref{Ultraproducto} we can consider an ultraproduct
${L_i^X}_{/{\cal U}}$ of the kind considered after Definition
\ref{LA}. We define $H(L_i) = {L_i^X}_{/{\cal U}}$. It is clear
that $H(L_i)$ is a linearly ordered ${\cal A}$-algebra. If we
take the class $e_i = [1_X]$ then $e_i$ is a unity in $H(L_i)$
since for every natural number $n$, $0 < e_i^n$ iff \hspace{0.1
cm} $\{x\in X : 0 < (1_X(x))^n\} \in {\cal U}$ and $\{x\in X : 0 <
(1_X(x))^n = x^n\} = X \in {\cal U}$.
\end{enumerate}

\noindent
 We can take the canonical embedding $j_i: L_i
\rightarrow H(L_i)$ and then for each $i\in I$ we can consider
$e_i$ as a unity lower bound of $L_i$ in $H(L_i)$. By Corollary
\ref{UNMTL}, $(e_i)_{i\in I}$ is a unity in $\prod_{i\in I}
H(L_i)$. Let  $j:\prod_{i\in I} L_i \rightarrow \prod_{i\in I}
H(L_i)$ be the monomorphism defined by $j((x_i)_{i\in I})=
(j_i(x_i))_{i\in I}$. Since $A$ is an absolute retract there
exists an epimorphism $\varphi:\prod_{i\in I} H(L_i)\rightarrow
A$ such that the following diagram commutes:

\begin{center}
\unitlength=1mm
\begin{picture}(90,20)(0,0)
\put(8,16){\vector(3,0){5}} \put(32,16){\vector(3,0){5}}

\put(52,10){\vector(0,-2){5}} \put(8,12){\vector(3,-1){40}}
\put(2,16){\makebox(0,0){$A$}}
\put(22,16){\makebox(0,0){$\prod_{i\in I} L_i$}}
\put(52,16){\makebox(0,0){$\prod_{i\in I} H(L_i)$}}

\put(52,-2){\makebox(0,0){$A$}}

\put(2,20){\makebox(17,0){$f$}} \put(26,20){\makebox(17,0){$j$}}

\put(28,9){\makebox(13,0){$\equiv$}}

\put(58,8){\makebox(-6,0){$\varphi$}}
\put(26,0){\makebox(-4,3){$1_A$}}

\end{picture}
\end{center}

\noindent
\\ Let $e = \varphi((e_i)_{i\in I})$. It is clear that $e$ is a
unity in $A$ since $\varphi$ is an homomorphism. If $u$ is a
unity in $A$ then $(e_i)_{i\in I}\leq jf(u)$ and by commutativity
of the above diagram, $e = \varphi((e_i)_{i\in I})\leq \varphi
jf(u)= u$. Thus $e = min\{u\in A : u \hspace{0.1 cm} is
\hspace{0.1 cm} unity\}$ resulting in $Rad(A) =  [e)$.

\qed
\end{proof}

\section  {Injectives in WNM-algebras and ${\cal MTL}$}

\begin{definition}
{\rm A {\it WNM-algebra} (weak nilpotent minimum) {\rm
\cite{GE1}} is an MTL-algebra satisfying the equation
$$\neg (x\odot y) \lor ((x\land y )\rightarrow (x\odot y)) = 1 . \leqno(W) $$}
\end{definition}

\noindent The variety of WNM-algebras is noted by ${\cal WNM}$.

\begin{theo}\label{SIMPLEWNM}
The following conditions are equivalent:

\begin{enumerate}
\item
$I$ is a simple WNM-algebra.

\item
$I$ has a coatom $u$ and its operations are given by

$$
x\odot y = \cases {0, & if $x, y < 1$\cr x, & if $y = 1$\cr y, &
if $x = 1$\cr}
$$

$$
x\rightarrow y = \cases {1, & if $x \leq y$ \cr y, & if $x = 1$\cr
u, & if $y < x < 1 . $\cr }
$$

\end{enumerate}
\end{theo}

\begin{proof}
$\Rightarrow$). For $Card(I)=2$ this result is trivial. If
$Card(I)>2$ then we only need to prove the following steps:

\begin{enumerate}

\item[a)]
 \textit{If $x, y < 1$ in $I$ then $x\odot y = 0$}: Since $I$ is
 simple, equation (W) implies that $x^2 = 0$ for each $x \in I
 \setminus \{1\}$. Hence if $x \leq y < 1$, then $x \odot y \leq y
 \odot y = 0$.
\item[b)]
\textit{$I$ has a coatom}: Let $0< x < 1$. We have that $\neg x <
1$ and, since $I$ is simple, we also have $\neg \neg x < 1$. Then
by a) it follows that $ \neg x \leq \neg \neg x \leq \neg \neg
\neg x = \neg x$, i. e., $\neg x = \neg \neg x$. If $0 < x, y <
1$, again by a) we have  $\neg x \odot \neg y = 0$. Thus $\neg x
\leq \neg \neg y = \neg y$. By interchanging $x$ and $y$ we
obtain the equality $\neg x = \neg y$. Now it is clear that  if
$0 < x <1$, then $u = \neg x$ is the coatom in $I$.

\item[c)]\textit{
If $y < x < 1$ then $x\rightarrow y = u$}: Since $x\rightarrow y =
\bigvee\{t\in I: t\odot x \leq y\}$, this supremum cannot be $1$
because $y < x $. Thus, in view of item a), $x\rightarrow y$ is
the coatom $u$.

\end{enumerate}
\noindent $\Leftarrow$) Immediate.

\qed
\end{proof}

\begin{example}
{\rm We can build simple WNM-algebras having arbitrary
cardinality if we consider an ordinal  $\gamma = Suc \hspace{0.05
cm} (Suc \hspace{0.05 cm}(\alpha))$ with the structure given by
Proposition \ref{SIMPLEWNM}, taking $Suc(\alpha)$ as coatom.
These algebras will be called {\it ordinal algebras}.}
\end{example}

\begin{prop}
${\cal WNM}$ and ${\cal MTL}$ have only trivial injectives.
\end{prop}

\begin{proof}
Follows from Proposition \ref{Injective simple} since these
varieties contain all ordinal algebras. \qed
\end{proof}

\section  {Injectives in SMTL-algebras}

\begin{definition}
{\rm An {\it SMTL-algebra}  {\rm \cite{GE2}} is a MTL-algebra
satisfying equation {\rm ($S$)}. The variety of SMTL-algebras is
denoted by ${\cal SMTL}$}.
\end{definition}

\begin{prop}
The only injectives in ${\cal SMTL}$ are complete boolean
algebras.
\end{prop}

\begin{proof}
Follows from Corollary \ref{SRLINJ} and  Theorem \ref{TEST1}
since the $test_d$-algebra $H_4$ belongs to ${\cal SMTL}$. \qed
\end{proof}

\section  {Injectives in $\Pi$SMTL-algebras}

\begin{definition}
{\rm A  {\it $\Pi$SMTL-algebra}  {\rm \cite{GE1}} is a
SMTL-algebra satisfying the equation:$$(\neg \neg z \odot
((x\odot z) \rightarrow (y\odot z)))\rightarrow (x\rightarrow y)
= 1. \leqno(\Pi)$$}
\end{definition}

\noindent The variety of $\Pi$SMTL-algebras is denoted by
$\Pi{\cal SMTL}$.

\begin{prop}\label{DENSIDP}
Let $A$ be an $\Pi$SMTL-algebra. Then $1$ is the only idempotent
dense element in $A$.
\end{prop}

\begin{proof}
By equation $\Pi$ it is easy to prove that, for each dense
element $\epsilon$, if $\epsilon \odot x = \epsilon \odot y$ then
$x = y$. Thus if  $\epsilon$ is an idempotent dense then
$\epsilon \odot 1 = \epsilon \odot \epsilon$ and $\epsilon = 1$.
\end{proof}

\begin{theo}\label{INJPI}
Let ${\cal A}$ be a subvariety of $\Pi{\cal SMTL}$. Then the
injectives in ${\cal A}$ are exactly the complete boolean
algebras.
\end{theo}

\begin{proof}
Follows from Corollary \ref{SRLINJ}, Theorem \ref{H3} and
Proposition \ref{DENSIDP}.
\end{proof}

\section  {Injectives in BL, MV, PL, and in Linear Heyting algebras}
\begin{definition}
{\rm A {\it BL-algebra}  {\rm \cite{HAJ}} is an MTL-algebra
satisfying the equation $$x\odot(x\rightarrow y) = x\land y
\leqno(B).$$}
\end{definition}

\noindent We denote by ${\cal BL}$ the variety of BL-algebras.
Important subvarieties of ${\cal BL}$ are the variety ${\cal MV}$
of multi-valued logic algebras (MV-algebras for short),
characterized by the equation  $\neg \neg x = x$ \cite{CDM,HAJ},
the variety ${\cal PL}$ of product logic algebras  (PL-algebras
for short), characterized by the equations ($\Pi$) plus ($S$)
\cite{HAJ,CT}, and the variety ${\cal HL}$ of linear Heyting
algebras, characterized by the equation $x\odot y = x\land y$
(also known as  G\"odel algebras \cite{HAJ}).

\begin{rema}\label{MV}
{\rm It is well known that ${\cal MV}$ is generated by the
MV-algebra $R_{[0,1]}= \langle [0,1], \odot, \rightarrow, \land,
\lor, 0, 1 \rangle$ such that $[0,1]$ is the real unit segment,
$\land$, $\lor$ are the natural meet and join on $[0,1]$ and
$\odot$ and $\rightarrow$ are defined as follows: $x\odot y:=
max(0,x+y-1)$, $x\rightarrow y:= min(1,1-x+y)$. $R_{[0,1]}$  is
the maximum simple algebra in ${\cal MV}$ {\rm (see \cite[Theorem
3.5.1]{CDM})}. Moreover $R_{[0,1]}$ is a rigid algebra {\rm (see
\cite[Corollary 7.2.6]{CDM})}, hence self-injective. Injective
MV-algebras were characterized in {\rm \cite[Corollary
2.11]{GLUS})} as the retracts of  powers of $R_{[0,1]}$.}
\end{rema}

\begin{prop}
If ${\cal A}$ is a subvariety of ${\cal PL}$, then the only
injectives of ${\cal A}$ are the complete boolean algebras.
\end{prop}

\begin{proof}
Follows from Theorem \ref{INJPI} since  ${\cal PL}$ is a
subvariety of $\Pi{\cal SMTL}$. \qed
\end{proof}

\begin{prop}
The only injectives in ${\cal HL}$ are the complete boolean
algebras.
\end{prop}

\begin{proof}
Follows from Corollary \ref{SRLINJ} and  Theorem \ref{TEST1}
since the algebra $test_d$ $H_4$ lies in ${\cal SMTL}$. \qed
\end{proof}

\begin{prop}\label{SIMBL}
${\cal BL}$ is a radical-dense variety.
\end{prop}

\begin{proof}
See {\rm \cite[Theorem 1.7 and Remark 1.9]{CT1}}. \qed
\end{proof}

\begin{prop}
Injectives in ${\cal BL}$ are exactly the retracts of  powers of
the MV-algebra $R_{[0,1]}$ .
\end{prop}

\begin{proof}
By  Remark \ref{MV} and Propositions \ref{SIMBL} and \ref{Simple
injective}, retracts of a power of the $R_{[0,1]}$ are injectives
in ${\cal BL}$. Thus by Theorem \ref{TEST1}, they are the only
possible injectives since $H_4$ lies in ${\cal BL}$. \qed
\end{proof}

\section  {Injectives in IMTL-algebras}

\begin{definition}
{\rm An involutive MTL-algebra (or {\it IMTL-algebra}) {\rm
\cite{GE1}} is a MTL-algebra satisfying the equation
$$\neg \neg x = x.  \leqno(I)$$}
\end{definition}

\noindent The variety of IMTL-algebras is noted by ${\cal IMTL}$. \\

An interesting IMTL-algebra, whose role is analogous to $H_3$ in
the radical-dense varieties, is the four element chain $I_4$
defined as follows:

\begin{center}
\unitlength=1mm
\begin{picture}(60,20)(0,0)

\put(-26,19){\makebox(0,0){$\odot$}} \put(-29,16){\line(3,0){38}}
\put(-23,20){\line(0,-2){30}}

\put(-28,19){\makebox(17,0){$1$}}
\put(-20,19){\makebox(17,0){$a$}}
\put(-11,19){\makebox(17,0){$b$}} \put(-3,19){\makebox(17,0){$0$}}

\put(-24,13){\makebox(-5,0){$1$}} \put(-24,6){\makebox(-5,0){$a$}}
\put(-24,-1){\makebox(-5,0){$b$}}
\put(-24,-8){\makebox(-5,0){$0$}}

\put(-17,13){\makebox(-5,0){$1$}} \put(-17,6){\makebox(-5,0){$a$}}
\put(-17,-1){\makebox(-5,0){$b$}}
\put(-17,-8){\makebox(-5,0){$0$}}

\put(-9,13){\makebox(-5,0){$a$}} \put(-9,6){\makebox(-5,0){$a$}}
\put(-9,-1){\makebox(-5,0){$0$}} \put(-9,-8){\makebox(-5,0){$0$}}

\put(0,13){\makebox(-5,0){$b$}} \put(0,6){\makebox(-5,0){$0$}}
\put(0,-1){\makebox(-5,0){$0$}} \put(0,-8){\makebox(-5,0){$0$}}

\put(8,13){\makebox(-5,0){$0$}} \put(8,6){\makebox(-5,0){$0$}}
\put(8,-1){\makebox(-5,0){$0$}} \put(8,-8){\makebox(-5,0){$0$}}

\put(21,19){\makebox(0,0){$\rightarrow$}}
\put(19,16){\line(3,0){38}} \put(25,20){\line(0,-2){30}}

\put(21,19){\makebox(17,0){$1$}} \put(28,19){\makebox(17,0){$a$}}
\put(36,19){\makebox(17,0){$b$}} \put(44,19){\makebox(17,0){$0$}}

\put(24,13){\makebox(-5,0){$1$}} \put(24,6){\makebox(-5,0){$a$}}
\put(24,-1){\makebox(-5,0){$b$}} \put(24,-8){\makebox(-5,0){$0$}}

\put(31,13){\makebox(-5,0){$1$}} \put(31,6){\makebox(-5,0){$1$}}
\put(31,-1){\makebox(-5,0){$1$}} \put(31,-8){\makebox(-5,0){$1$}}

\put(39,13){\makebox(-5,0){$a$}} \put(39,6){\makebox(-5,0){$1$}}
\put(39,-1){\makebox(-5,0){$1$}} \put(39,-8){\makebox(-5,0){$1$}}

\put(47,13){\makebox(-5,0){$b$}} \put(47,6){\makebox(-5,0){$b$}}
\put(47,-1){\makebox(-5,0){$1$}} \put(47,-8){\makebox(-5,0){$1$}}

\put(55,13){\makebox(-5,0){$0$}} \put(55,6){\makebox(-5,0){$b$}}
\put(55,-1){\makebox(-5,0){$a$}} \put(55,-8){\makebox(-5,0){$1$}}

\put(72,20){\line(0,-2){30}} \put(72,20){\circle*{1.5}}
\put(72,10){\circle*{1.5}} \put(72,0){\circle*{1.5}}
\put(72,-10){\circle*{1.5}}

\put(78,20){\makebox(-5,0){$1$}} \put(78,10){\makebox(-5,0){$a$}}
\put(83,-0){\makebox(-5,0){$b= \neg a$}}
\put(78,-10){\makebox(-5,0){$0$}}

\end{picture}
\end{center}

\vspace{1cm}

\begin{theo}\label{I4}
Let ${\cal A}$ be a a subvariety of ${\cal IMTL}$. If $A$ is a
non-semisimple absolute retract in ${\cal A}$, then $A$ has a
principal unity $\epsilon$ and  $\{0, \neg \epsilon,
\epsilon,1\}$ is a subalgebra of $A$ which is isomorphic to $I_4$.
\end{theo}

\begin{proof}
Follows from Theorem \ref{MINUNITY}. \qed
\end{proof}

\begin{definition}
{\rm Let ${\cal A}$ be a subvariety of ${\cal IMTL}$. An algebra
$T$ is called {\it $test_I$-algebra} iff, it has a subalgebra
$\{0,\neg \epsilon, \epsilon, 1\}$ isomorphic to $I_4$ and there
exists $t \in Rad(T)$ such that $t < \epsilon$}.
\end{definition}

\begin{theo}\label{TEST2}
Let ${\cal A}$ be a subvariety of ${\cal IMTL}$. If ${\cal A}$
has a nontrivial injective and contains a $test_I$-algebra, then
injectives are semisimple.
\end{theo}

\begin{proof}
Let T be a $test_I$-algebra and  $t\in Rad(T_i)$ such that $t <
\epsilon$. We can consider  a subdirect embedding $f:T \rightarrow
\prod_{j\in J} H_j$ such that $L_j$ is linearly ordered. Let $x_j
= \pi_j f(x)$ for each $x\in T$ and $\pi_j$ the {\it j}{\rm
th}-projection. Since $t < \epsilon$, exists $s\in J$ such that
$\neg \epsilon_s < \neg t_s < t_s < \epsilon_s$ and by Corollary
\ref{UNMTL} , $t_s$ and $\epsilon_s$ are unities in the chain
$H_s$ with $\epsilon_s$ idempotent. Note that $H_s$ is also a
$test_I$-algebra. To see that $\epsilon_s \rightarrow t_s \leq
\epsilon$, observe first that $0 < \epsilon_s \odot \neg t_s$
since, if $\epsilon_s \odot \neg t_s = 0$ then $\epsilon_s \leq
\neg \neg t_s = t_s$ which is a contradiction. Consequently,
$\neg \epsilon_s \leq \epsilon_s \odot \neg t_s$ since, if $
\epsilon_s \odot \neg t_s  \leq \neg \epsilon_s$ then $\epsilon_s
\odot \neg t_s =  (\epsilon_s)^2 \odot \neg t_s  \leq  \neg
\epsilon \odot \epsilon = 0$. Thus we can conclude that
$\epsilon_s \rightarrow t_s = \neg(\epsilon_s \odot \neg t_s)
\leq \neg \neg \epsilon_s = \epsilon_s$. Suppose that there
exists a non-semisimple injective $A$ in ${\cal A}$. Then by
Theorem \ref{I4}, let $\alpha: I_4 \rightarrow A$ be a
monomorphism such that $\alpha(a)$ is the principal unity in $A$.
Let $i:I_4\rightarrow H_s$ be the monomorphism such that $i(a)=
\epsilon_s$. Since $A$ is injective, there exists a homomorphism
$\varphi:H_s\rightarrow A$ such that the following diagram
commutes:

\begin{center}
\unitlength=1mm
\begin{picture}(20,20)(0,0)
\put(8,16){\vector(3,0){5}} \put(2,10){\vector(0,-2){5}}
\put(8,2){\vector(1,1){7}}

\put(2,10){\makebox(13,0){$\equiv$}}

\put(2,16){\makebox(0,0){$I_4$}} \put(20,16){\makebox(0,0){$A$}}
\put(2,0){\makebox(0,0){$H_s$}}
\put(2,20){\makebox(17,0){$\alpha$}}
\put(2,8){\makebox(-6,0){$i$}}
\put(16,0){\makebox(-4,3){$\varphi$}}
\end{picture}
\end{center}

\noindent Since $\alpha(a)$ is the principal unity in $A$ and
$t_s\leq \epsilon_s$ then, by commutativity, $\varphi(\epsilon_s)
= \varphi(t_s)= \alpha(a)$. Thus $\varphi(\epsilon_s \rightarrow
t_s) = 1$, which is a contradiction since $\varphi(\epsilon_s
\rightarrow t_s) \leq \varphi(\epsilon_s) = \alpha(a) < 1$. Hence
${\cal A}$ has only semisimple injectives. \qed
\end{proof}

\begin{prop}\label{NOINJIMTL}
${\cal IMTL}$ has only trivial injectives.
\end{prop}

\begin{proof}
Suppose that there exists nontrivial injectives in ${\cal IMTL}$.
By Theorem \ref{Injective simple} there is a simple maximum
algebra $I$ in ${\cal IMTL}$. We consider the six elements $IMTL$
chain $I_6$ defined as follows:

\begin{center}
\unitlength=1mm
\begin{picture}(60,20)(0,0)

\put(-29,19){\makebox(0,0){$\odot$}} \put(-32,16){\line(3,0){53}}
\put(-26,20){\line(0,-2){46}}

\put(-31,19){\makebox(17,0){$1$}}
\put(-23,19){\makebox(17,0){$a_1$}}
\put(-14,19){\makebox(17,0){$t$}}
\put(-6,19){\makebox(17,0){$a_2$}}
\put(2,19){\makebox(17,0){$a_3$}} \put(10,19){\makebox(17,0){$0$}}

\put(-27,13){\makebox(-5,0){$1$}}
\put(-27,6){\makebox(-5,0){$a_1$}}
\put(-27,-1){\makebox(-5,0){$t$}}
\put(-27,-8){\makebox(-5,0){$a_2$}}
\put(-27,-15){\makebox(-5,0){$a_3$}}
\put(-27,-22){\makebox(-5,0){$0$}}

\put(-20,13){\makebox(-5,0){$1$}}
\put(-20,6){\makebox(-5,0){$a_1$}}
\put(-20,-1){\makebox(-5,0){$t$}}
\put(-20,-8){\makebox(-5,0){$a_2$}}
\put(-20,-15){\makebox(-5,0){$a_3$}}
\put(-20,-22){\makebox(-5,0){$0$}}

\put(-12,13){\makebox(-5,0){$a_1$}}
\put(-12,6){\makebox(-5,0){$a_2$}}
\put(-12,-1){\makebox(-5,0){$a_3$}}
\put(-12,-8){\makebox(-5,0){$a_3$}}
\put(-12,-15){\makebox(-5,0){$0$}}
\put(-12,-22){\makebox(-5,0){$0$}}

\put(-3,13){\makebox(-5,0){$t$}} \put(-3,6){\makebox(-5,0){$a_3$}}
\put(-3,-1){\makebox(-5,0){$a_3$}}
\put(-3,-8){\makebox(-5,0){$0$}} \put(-3,-15){\makebox(-5,0){$0$}}
\put(-3,-22){\makebox(-5,0){$0$}}

\put(5,13){\makebox(-5,0){$a_2$}} \put(5,6){\makebox(-5,0){$a_3$}}
\put(5,-1){\makebox(-5,0){$0$}} \put(5,-8){\makebox(-5,0){$0$}}
\put(5,-15){\makebox(-5,0){$0$}} \put(5,-22){\makebox(-5,0){$0$}}

\put(13,13){\makebox(-5,0){$a_3$}} \put(13,6){\makebox(-5,0){$0$}}
\put(13,-1){\makebox(-5,0){$0$}} \put(13,-8){\makebox(-5,0){$0$}}
\put(13,-15){\makebox(-5,0){$0$}}
\put(13,-22){\makebox(-5,0){$0$}}

\put(21,13){\makebox(-5,0){$0$}} \put(21,6){\makebox(-5,0){$0$}}
\put(21,-1){\makebox(-5,0){$0$}} \put(21,-8){\makebox(-5,0){$0$}}
\put(21,-15){\makebox(-5,0){$0$}}
\put(21,-22){\makebox(-5,0){$0$}}

\put(30,19){\makebox(0,0){$\rightarrow$}}
\put(28,16){\line(3,0){53}} \put(34,20){\line(0,-2){46}}

\put(29,19){\makebox(17,0){$1$}}
\put(38,19){\makebox(17,0){$a_1$}}
\put(46,19){\makebox(17,0){$t$}}
\put(54,19){\makebox(17,0){$a_2$}}
\put(62,19){\makebox(17,0){$a_3$}}
\put(70,19){\makebox(17,0){$0$}}

\put(32,13){\makebox(-5,0){$1$}} \put(32,6){\makebox(-5,0){$a_1$}}
\put(32,-1){\makebox(-5,0){$t$}}
\put(32,-8){\makebox(-5,0){$a_2$}}
\put(32,-15){\makebox(-5,0){$a_3$}}
\put(32,-22){\makebox(-5,0){$0$}}

\put(40,13){\makebox(-5,0){$1$}} \put(40,6){\makebox(-5,0){$1$}}
\put(40,-1){\makebox(-5,0){$1$}} \put(40,-8){\makebox(-5,0){$1$}}
\put(40,-15){\makebox(-5,0){$1$}}
\put(40,-22){\makebox(-5,0){$1$}}

\put(49,13){\makebox(-5,0){$a_1$}} \put(49,6){\makebox(-5,0){$1$}}
\put(49,-1){\makebox(-5,0){$1$}} \put(49,-8){\makebox(-5,0){$1$}}
\put(49,-15){\makebox(-5,0){$1$}}
\put(49,-22){\makebox(-5,0){$1$}}

\put(57,13){\makebox(-5,0){$t$}} \put(57,6){\makebox(-5,0){$a_1$}}
\put(57,-1){\makebox(-5,0){$1$}} \put(57,-8){\makebox(-5,0){$1$}}
\put(57,-15){\makebox(-5,0){$1$}}
\put(57,-22){\makebox(-5,0){$1$}}

\put(65,13){\makebox(-5,0){$a_2$}}
\put(65,6){\makebox(-5,0){$a_1$}}
\put(65,-1){\makebox(-5,0){$a_1$}}
\put(65,-8){\makebox(-5,0){$1$}} \put(65,-15){\makebox(-5,0){$1$}}
\put(65,-22){\makebox(-5,0){$1$}}

\put(73,13){\makebox(-5,0){$a_3$}} \put(73,6){\makebox(-5,0){$t$}}
\put(73,-1){\makebox(-5,0){$a_1$}}
\put(73,-8){\makebox(-5,0){$a_1$}}
\put(73,-15){\makebox(-5,0){$1$}}
\put(73,-22){\makebox(-5,0){$1$}}

\put(81,13){\makebox(-5,0){$0$}} \put(81,6){\makebox(-5,0){$a_3$}}
\put(81,-1){\makebox(-5,0){$a_2$}}
\put(81,-8){\makebox(-5,0){$t$}}
\put(81,-15){\makebox(-5,0){$a_1$}}
\put(81,-22){\makebox(-5,0){$1$}}

\put(88,20){\line(0,-2){50}} \put(88,20){\circle*{1.5}}
\put(88,10){\circle*{1.5}} \put(88,0){\circle*{1.5}}
\put(88,-10){\circle*{1.5}} \put(88,-20){\circle*{1.5}}
\put(88,-30){\circle*{1.5}}

\put(94,20){\makebox(-5,0){$1$}}
\put(94,10){\makebox(-5,0){$a_1$}}
\put(94,-0){\makebox(-5,0){$t$}}
\put(94,-10){\makebox(-5,0){$a_2$}}
\put(94,-20){\makebox(-5,0){$a_3$}}
\put(94,-30){\makebox(-5,0){$0$}}

\end{picture}
\end{center}

\vspace{3cm}

\noindent Since $I$ is simple maximum we can consider $I_6$ and
$R_{[0,1]}$ as subalgebras of $I$. In view of this and using the
nilpotence order we have that $1/2 < t < 3/4$ since $I$ is a
chain. Therefore we can consider $u = \bigvee_{R_{[0,1]}} \{ x\in
R_{[0,1]}:x < t\}$ and $v = \bigwedge _{R_{[0,1]}}\{x\in
R_{[0,1]}:x > t\}$ and it is clear that $u,v \in R_{[0,1]}$ since
$R_{[0,1]}$ is a complete algebra. Thus $u < t < v$. This
contradicts the fact that the order of $R_{[0,1]}$ is dense.
Consequently ${\cal IMTL}$ has only trivial injectives. \qed
\end{proof}

\section  {Injectives in NM-algebras}

\begin{definition}
{\rm A nilpotent minimum algebra (or {\it NM-algebra}) {\rm
\cite{GE1}} is an IMTL-algebra satisfying the equation {\rm
($W$)}}.
\end{definition}

\noindent The variety of NM-algebras is noted by ${\cal NM}$.  As
an example we consider $N_{[0,1]}= \langle [0,1], \odot,
\rightarrow, \land, \lor, 0, 1 \rangle$ such that $[0,1]$ is the
real unit segment, $\land$, $\lor$ are the natural meet and join
on $[0,1]$ and $\odot$ and $\rightarrow$ are defined as follows:

$$
x\odot y = \cases {x\land y, & if $1 < x + y $\cr 0, & otherwise,
\cr}
$$

$$
x\rightarrow y = \cases {1, & if $x \leq y$ \cr max (y, 1-x) &
otherwise .\cr }
$$

\vspace{0.5cm}

\noindent Note that $\{0, \frac{1}{2}, 1\}$ is the universe of a
subalgebra of $N_{[0,1]}$, that we denote by $\L_3$. The
subvariety of ${\cal NM}$ generated by $\L_3$ coincides with the
variety ${\cal L}_3$ of three-valued \L ukasiewicz algebras (see
\cite{Mo,Cig}).

\begin{prop}\label{SIMPLENM}
$\L_3$ is the maximum simple algebra in ${\cal NM}$, and  it is
self-injective.
\end{prop}

\begin{proof}
Let $I$ be a simple algebra such that $Card(I)>2$. By Theorem
\ref{SIMPLEWNM} $I$ has a coatom $u$ satisfying $\neg x = u$ for
each $0 < x <1$. Thus $x = \neg \neg x = \neg u = u$ for each $0
< x <1$. Consequently $Card(I)= 3$ and $I = \L_3$. \qed
\end{proof}

\begin{coro}\label{LUK}
${\cal S}em({\cal NM}) = {\cal L}_3$. \qed
\end{coro}

\begin{prop}
Injectives in ${\cal NM}$ coincide with complete Post algebras of
order $3$.
\end{prop}

\begin{proof}
By Proposition \ref{SIMPLEWNM}, Theorem  \ref{Simple injective}
and Theorem  \ref{TEST2} injectives in ${\cal NM}$ are semisimple
since $N_{[0,1]}$ is an algebra $Test_I$. Thus by Proposition
\ref{LUK} and {\rm \cite{Mo}, \cite[Theorem 3.7]{Cig}}, complete
Post algebras of order $3$ are the injectives in ${\cal NM}$. \qed
\end{proof}

\section  {Injective bounded hoops}

\begin{definition}
{\rm A {\it hoop} {\rm \cite{BlokFer}} is an algebra $ \langle A,
\odot, \rightarrow, 1 \rangle$ of type $ \langle 2, 2, 0 \rangle$
satisfying the following axioms:

\begin{enumerate}
\item
$\langle A,\odot,1 \rangle$ is an abelian monoid,

\item
$x\rightarrow x = 1$,

\item
$(x\rightarrow y)\odot x = (y\rightarrow x)\odot y$,

\item
$x\rightarrow (y\rightarrow z) = (x\odot y)\rightarrow z$.

\end{enumerate}
}
\end{definition}

\noindent The variety of hoops is noted ${\cal HO}$. Every hoop
is a meet semilattice, where the meet operation is given by
$x\land y = x\odot (x\rightarrow y)$. Let $A$ be a hoop. If $A$
has smallest element $0$, we can define an unary operation $\neg$
by $\neg x = x\rightarrow 0$.  A subset $F$ of $A$ is a {\bf
filter} iff $1\in F$ and  $F$ is closed under $\odot$. As in
residuated lattices, filters and congruences can be identified
{\rm \cite {BlokFer}}.

\begin{definition}
{\rm A {\it Wajsberg hoop} {\rm \cite {BlokFer}} is a hoop that
satisfies the following equation $$(x\rightarrow y) \rightarrow y
= (y\rightarrow x)\rightarrow x.  \leqno(T)$$}
\end{definition}

\noindent Each Wajsberg hoop is a lattice, in which the join
operation is given by $x\lor y = (x\rightarrow y) \rightarrow y$.

\begin{definition}
{\rm A {\it bounded hoop} is an algebra $ \langle A, \odot,
\rightarrow, 0, 1 \rangle$ of type $ \langle 2, 2, 0, 0 \rangle$
such that:

\begin{enumerate}

\item
$ \langle A, \odot, \rightarrow, 1 \rangle$ is a hoop

\item
$0\rightarrow x = 1.$

\end{enumerate}
}
\end{definition}

The variety of bounded hoop is noted by ${\cal BH}_0$. Observe that since $0$ is in the clone of hoop operation, we require that for each morphism $f$, $f(0)=0$. In the same way as in the case of residuated lattices, for each bounded hoop $A$, we can consider $Ds(A)$ the set of dense elements of $A$, and this is an implicative filter of $A$. \\

\begin{prop}\label{SIMPLEHOOP}
A bounded simple hoop is a simple MV-algebra.
\end{prop}

\begin{proof}
Let $I$ be a simple hoop. Then by {\rm \cite [Corollary
2.3]{BlokFer}} it is a totally ordered Wajsberg hoop. If $0$ is
the smallest element in $I$ then by the equation (T), $ \neg \neg
x = (x\rightarrow 0)\rightarrow 0 = (0\rightarrow x) \rightarrow
x = 1\rightarrow x = x$. Hence it is an MV-algebra. Since the
MV-congruences are in correspondence with implicative filters,
$I$ is a simple MV-algebra. \qed
\end{proof}

\begin{table}[h] \begin{center}{\scriptsize
\begin{tabular}{|l|l|l|}\hline
\multicolumn{1}{|c|}{Variety} & \multicolumn{1}{|c|}{Equations} &
\multicolumn{1}{|c|}{Injectives} \\
\hline\hline ${\cal RL}$ &
 &
Trivial \\ \hline

${\cal DRL}$ & ${\cal RL} + x\land (y\lor z) = (x\land y) \lor
(x\land z)$ & Trivial  \\ \hline

${\cal GM}$ & ${\cal RL} + \neg \neg x = x$ & Trivial  \\ \hline

${\cal DGM}$ & ${\cal GM} + x\land (y\lor z) = (x\land y) \lor
(x\land z)$ & Trivial  \\ \hline

${\cal MTL}$ & ${\cal RL} + (x\rightarrow y) \lor (y\rightarrow
x) = 1$ & Trivial  \\ \hline

${\cal WNM}$ & ${\cal MTL} + \neg(x\odot y)\lor ((x\land y)
\rightarrow (x\odot y ))=1 $ & Trivial  \\ \hline

${\cal IMTL}$ & ${\cal MTL} + \neg \neg x = x$ & Trivial  \\
\hline

${\cal BL}$ & ${\cal MTL} + x\land y = x \odot (x\rightarrow y)$ &
Retracts of  powers of $R_{[0,1]}$ \\ \hline

${\cal MV}$ & ${\cal BL} + \neg \neg x = x$ & Retracts of  powers
of $R_{[0,1]}$ \\ \hline

${\cal BH}_0$ & $\lor$-free subreduct of ${\cal RL} + x\land y =
x \odot (x\rightarrow y)$  & Retracts of  powers of $R_{[0,1]}$ \\
\hline

${\cal SRL}$ & ${\cal RL} + x\land \neg x = 0$ & Complete boolean
algebras  \\ \hline

${\cal SMTL}$ & ${\cal MTL} + x\land \neg x = 0$ & Complete
boolean algebras  \\ \hline

$\Pi{\cal SMTL}$ & ${\cal SMTL} + \neg \neg z \odot ((x\odot
z)\rightarrow (y\odot z))\leq (x\rightarrow y )$ & Complete
boolean algebras  \\ \hline

${\cal PL}$ & $\Pi{\cal SMTL} + x\land y = x \odot (x\rightarrow
y)$ & Complete boolean algebras  \\ \hline

${\cal HL}$ & ${\cal BL} + x\land y = x \odot y$ & Complete
boolean algebras  \\ \hline

${\cal NM}$ & ${\cal WNM} + \neg \neg x = x$ & Complete Post
algebras of order 3  \\ \hline

\end{tabular}}
\caption {Injectives in Varieties of Residuated Algebras}
\end{center}
\end{table}

\begin{prop}\label{SIMPLEHOM}
Let $I,J$ be simple hoops with smallest elements $0_I, 0_J$
respectively. If $\varphi: I \rightarrow J$ is a hoop
homomorphism then $\varphi$ is also an MV-homomorphism, i.e.,
$\varphi(0_I)=0_J$.
\end{prop}

\begin{proof}
Suppose that $\varphi(0_I) = a$. Since $J$ is simple, there
exists a natural number $n$ such that $a^n = 0_J$. Thus we have,
$\varphi(0_I) = \varphi(0_I^n) = (\varphi(0_I))^n = a^n = 0_J$.
\qed
\end{proof}

\noindent The following two results are obtained in the same way
as Theorems \ref{H3} and \ref{TEST1} respectively.

\begin{theo}\label{HH3}
Let ${\cal A}$ be a subvariety of ${\cal BH}_0$. If $A$ is a
non-semisimple absolute retract in ${\cal A}$, then $Ds(A)$ has a
least element $\epsilon$ i.e, $Ds(A) = [\epsilon)$ and
$\{0,\epsilon,1\}$ is a subalgebra of $A$ isomorphic to the three
element Heyting algebra $H_3$. \qed
\end{theo}

\begin{theo}\label{TEST3}
Let ${\cal A}$ be a subvariety of ${\cal BH}_0$. If ${\cal A}$
has a nontrivial injectives and contains the Heyting algebra
$H_4$ then injectives are semisimple. \qed
\end{theo}

\begin{coro}
Injectives in ${\cal BH}_0$ are exactly the retracts of  powers of
the MV-algebra $R_{[0,1]}$.
\end{coro}

\begin{proof}
By Proposition \ref{SIMPLEHOOP}, semisimple bounded hoops are
MV-algebras. Therefore $R_{[0,1]}$ is the maximum simple algebra
and it is self injective by Proposition \ref{SIMPLEHOM}. Thus by
Theorem \ref{Simple injective} retracts of powers of the
MV-algebra $R_{[0,1]}$ are injectives in ${\cal BH}_0$. By
Theorem  \ref{TEST3} they are the only injectives, because $H_4$
lies in ${\cal BH}_0$. \qed
\end{proof}

{\small \noindent Departamento de Matem\'atica\\
Facultad de Ciencias Exactas y Naturales\\
Universidad de Buenos Aires\\
Ciudad Universitaria\\
1428 Buenos Aires - Argentina\\
e-mail: hfreytes@dm.uba.ar}

\end{document}